\documentclass[a4paper,12pt]{article}

\usepackage{amsmath}
\usepackage{amsthm}
\usepackage{amssymb}
\usepackage{latexsym}
\usepackage{graphicx}
\usepackage{stmaryrd}
\usepackage{mathrsfs}
\usepackage{enumerate}
\usepackage{cases}

\newcommand{\glalign}[2]{\lower.6ex\vbox{
\baselineskip\lineskip\ialign{$#1\hfil##\hfil$\crcr#2\crcr=\crcr}}}

\newcommand{\del}{\partial}

\renewcommand{\div}{\mbox{\rm div}\,}

\newcommand{\trans}{{}^\top}

\def\eqn#1$$#2$${\begin{equation}\label#1#2\end{equation}}

\numberwithin{equation}{section}

\newtheorem{defi}{Definition}[section]
\newtheorem{thm}[defi]{Theorem}

\newtheorem{prop}[defi]{Proposition}
\newtheorem{lem}[defi]{Lemma}
\newtheorem{rem}[defi]{Remark}

\def\eqn#1$$#2$${\begin{equation}\label#1#2\end{equation}}

\numberwithin{equation}{section}
\numberwithin{equation}{section}
\allowdisplaybreaks[4]

\begin{document}


\title{
\bf \large 
Time decay estimate with diffusion wave property and smoothing effect for solutions to the compressible Navier-Stokes-Korteweg 
system}
\author{
Takayuki KOBAYASHI and Kazuyuki TSUDA\\
 Osaka University, \\
 1-3, Machikaneyamacho, Toyonakashi, 560-8531, JAPAN \\
e-mail: jtsuda@sigmath.es.osaka-u.ac.jp}
\date{}
\maketitle
\begin{abstract} 
Time decay estimate of solutions to the compressible Navier-Stokes-Korteweg system is studied. Concerning the linearized problem, the decay estimate with diffusion wave property for an initial data is derived. As an application, the time decay estimate of solutions to the nonlinear problem is given.  In contrast to the compressible Navier-Stokes system, for linear  system 
regularities of the initial data are lower and independent of the order of derivative of solutions owing to smoothing effect from the Korteweg tensor.  Furthermore, for the nonlinear system diffusion wave property is obtained with  an initial data having lower regularity than that of study of the compressible Navier-Stokes system.  
\end{abstract}




\section{Introduction}
We study time decay estimate for solutions to the following compressible Navier-Stokes-Korteweg system 
in $\mathbb{R}^n$ $(n \geq 2)$:
\begin{numcases}
{}
\partial_{t}\rho +\div M=0,\nonumber\\
\partial_{t}M +\div \Big(\frac{M \otimes M}{\rho}\Big)+\nabla P(\rho)=\div \Big({\mathcal S}(\frac{M}{\rho})+{\mathcal K}(\rho)\Big),\label{CNSK}\\
\rho(x,0) =\rho_0, \ \ M(x,0)=M_0. \nonumber
\end{numcases}
\noindent Here $\rho=\rho(x,t)$ and $M=(M_{1}(x,t),\cdots,M_{n}(x,t))$ denote the unknown density and momentum respectively, at time $t\in \mathbb{R}^{+}$ and 
position $x\in\mathbb{R}^n$; $\rho_0=\rho_0(x)$ and $M_0=M_0(x)$ denote given initial datas; 
${\mathcal S}$ and ${\mathcal K}$ denote the viscous stress tensor and the Korteweg stress tensor that are given by 
\begin{eqnarray}
\left\{
\begin{array}{ll}
{\mathcal S}(\frac{M}{\rho}) =\Big(\mu' \div \frac{M}{\rho} \Big)\delta_{i,j} +2\mu d_{ij}\Big(\frac{M}{\rho}\Big),\\
{\mathcal K}(\rho) =\frac{\kappa}{2}(\Delta \rho^2 -|\nabla \rho|^2)\delta_{i,j}- \kappa \frac{\del \rho}{\del x_i}\frac{\del \rho}{\del x_j},
\end{array}
\right.\label{Korteweg tensor}
\end{eqnarray}
where $d_{ij}\Big(\frac{M}{\rho}\Big)=\frac{1}{2}\left(\frac{\del }{\del x_i}\Big(\frac{M}{\rho}\Big)_j+\frac{\del}{\del x_j}\Big(\frac{M}{\rho}\Big)_i\right)$; $\mu$ and $\mu'$ are the viscosity coefficients that are assumed to be constants satisfying 
$$
\mu>0, \ \ \ 
\frac{2}{n}\mu+\mu'\geq 0. 
$$   
$\kappa$ denotes the capillary constant that is assumed to be a positive constant. Note that if $\kappa=0$ in the Korteweg tensor, the usual compressible Navier-Stokes equation (the abbreviation is used by ``CNS'' below) is obtained; 
$P=P(\rho)$ is the pressure that is assumed to be a smooth function of $\rho$ satisfying  
$
P'(\rho_*)>0, 
$
where $\rho_*$ is a given positive constant and $(\rho_*,0)$ denotes a given constant state. 
We consider solutions to \eqref{CNSK} around the constant state. 

The system \eqref{CNSK} describes two phase flow
with phase transition between liquid and vapor in a compressible fluid as a diffuse
interface model.  
In the diffuse interface model, the phase boundary is regarded as a narrow transition layer and fluid state is described by change of the density.    
Hence it is enough to consider one set of equations and  a single spatial domain in contrast to the classical sharp interface model. It is well known that the phase field method use the idea of diffuse interface  effectively for numerical simulation. 

Concerning derivation of \eqref{CNSK}, 
Van der Waals \cite{Van der Waals} observes that a phase transition boundary can be regarded as a thin transition zone, i.e, diffuse interface caused by a steep gradient of the density.  Based on his idea, Korteweg \cite{Korteweg}  suggests the stress tensor including the term $\nabla\rho \otimes \nabla \rho$ of the Navier-Stokes equation. 
Then Dunn and Serrin \cite{Dunn and Serrin} generalize the Korteweg's work and provide the system $(\ref{CNSK})$ with (\ref{Korteweg tensor}).
In recent works, Heida and M\'{a}lek \cite{Heida and malek} 
derive \eqref{CNSK} by  the entropy production method in difference from \cite{Dunn and Serrin}.  
Freist\"{u}hler and Kotschote \cite{Kotschote} derive the Navier-Stokes-Allen-Cahn system and the Navier-Stokes-Chan-Hilliard system which describe two phase flow of mixture materials from some model of Korteweg type. Gorban and Karlin \cite{Gorban-Karlin} derive the Korteweg tensor from the Bolzmann equation.

We first study time decay estimate for solutions to linearized problem of \eqref{CNSK}. 
We shall show that the leading part of solutions consists of $\rho_*$ and a divergence free momentum field decaying in the same order as a $n$-dimensional heat kernel as $t$ goes to infinity in $L^p$ $(p\geq 2)$. 
Solutions are decomposed into low frequency part and high frequency part. 
We also shall show that  the low frequency part of density  decays faster than the 
$n$-dimensional heat kernel in $L^\infty$ norm and the decay order is obtained. On the other hand, we show that solutions may grow in $L^1$ norm as $t$ goes to infinity and the growth order is obtained. These properties are called as ``diffusion wave property" which occurs from 
terms in the Green matrix given by the convolution of the Green functions 
of the diffusion equation and the wave equation. 
The diffusion wave property is studied for CNS by Hoff and Zumbrun \cite{Hoff-Zumbrun1, Hoff-Zumbrun2}, and Kobayashi and Shibata \cite{Kobayashi-Shibata}. We also give $L^p$-$L^q$ $(1\leq q \leq 2 \leq p)$ estimate of solutions for the low frequency part.  Concerning the high frequency part, it is shown that solutions have exponential decay as $t\rightarrow \infty$ similarly to CNS.  It differs from CNS that smoothing effect of solutions appears in the estimate of the high frequency part.  

To show the decay estimate of solutions to the linearized problem, we use the Fourier transform method as in \cite{Kobayashi-Shibata} for CNS and Shibata \cite{Shibata-v} for the viscoelastic equation. In contrast to \cite{Kobayashi-Shibata}, due to the Korteweg tensor, the smoothing effect of the heat kernel appears in every components of the Green matrix of \eqref{CNSK}. Therefore we do not assume any regularity for the initial value in the estimate of the high frequency part. (See Theorem \ref{est-linear-high}  below. ) On the other hand, though the Korteweg tensor is added, the order of roots for characteristic equations of the linearized system 
coincides with that of CNS on the low frequency part in the Fourier space. Hence, the estimates of the low frequency part are derived similarly to the proofs of  \cite{Kobayashi-Shibata}. 

Furthermore, we derive time decay estimates for solutions to the nonlinear problem \eqref{CNSK}. Concerning \eqref{CNSK}, Danchin and Desjardins \cite{Danchin} show global existence of a solution around the motionless state $(\rho_*, 0)$ with a small initial data $u_0 \in (B^{\frac{n}{2}}_{2,1}\cap B^{\frac{n}{2}-1}_{2,1})\times B^{\frac{n}{2}-1}_{2,1}$, where $B^{\frac{n}{2}}_{2,1}$ denotes the usual homogeneous Besov space. 
Hattori and Li \cite{Hattori-Li-1,Hattori-Li-2} show the global existence of a $H^{N+1} \times H^{N}$ solution  with a small 
$u_0 \in H^{N+1}\times H^N$, where $N$ is an integer satisfying that $N\geq [n/2]+2$ and $[n/2]$ denotes the integer part of $n/2$. For three dimensional case Tan and R. Zhang, X. Zhang and Tan \cite{Tan-Zhang,Zhang-Tan} study the global existence for a small $u_0 \in H^{4}\times H^3$ with large time behavior of solutions.  Tan, Wang and Xu \cite{Tan-Wang-Xu} state the global existence of a solution which has lower regularity, that is, $C([0, \infty); H^{2} \times H^{1})$ class with a small 
$u_0 \in H^{2}\times H^1$. 
Wang and Tan \cite{Wang-Tan} show convergence rates of $L^p$ $(2\leq p)$ norms of the solution for $3$ dimensional case under a small initial value; 
Let the velocity field $v$ be defined by $v=M/\rho$. If $\|(\rho_0, v_0)\|_{H^{s+1}\times H^s}<<1$ $(s\geq 3)$, then we have that for $t>0$ 
\begin{eqnarray*}
&&\|(\rho(t)-\rho_*, v(t))\|_{L^p}\leq C(1+t)^{-\frac{3}{2}(1-\frac{1}{p})}  \ \ (2\leq p \leq 6),\\[1ex]
&&\|\nabla (\rho(t)-\rho_*, v(t))\|_{L^2}\leq C(1+t)^{-\frac{5}{4}},
\end{eqnarray*}
where $(\rho_0, v_0)$ denotes a given initial data and $H^s$ denotes the usual $L^2$ Sobolev space. If in addition $s\geq 4$ we have that 
$$
\|(\rho(t)-\rho_*, v(t))\|_{L^p}\leq C(1+t)^{-\frac{3}{2}(1-\frac{1}{p})} \ \ (2\leq p \leq \infty).
$$
Okita \cite{Okita} show that for $n$ dimensional case $L^2$ time decay rate of solutions to CNS around some stationary solution is obtained with a small $H^{s}$ initial data, where $s$ is an integer satisfying that $s\geq [n/2]+1$. 
We shall show that for \eqref{CNSK} in $n$ dimensional case nonlinear parts in the Duhamel formula have faster decay than linear parts and $L^\infty$ and $L^1$ estimates of solutions are stated as $t$ goes to infinity with the decay and growth rates. 
We also obtain $L^2$ time decay estimate of solutions with the decay rate similar to \cite{Wang-Tan} and CNS \cite{Okita}. 
It differs from \cite{Okita} and \cite{Wang-Tan}  that in our results the estimates reflect the diffusion wave property of the system. In both \cite{Okita} and \cite{Wang-Tan}  the diffusion wave property of solutions to the nonlinear problem is not studied.  In addition, by decomposition method for solution (cf. \cite{Okita} for CNS), the initial data is assumed in 
lower regularity class than that of \cite{Hoff-Zumbrun1, Wang-Tan}. Concerning $L^\infty$ estimate of solutions in the high frequency part, since that of linear problem has singularity at $t=0$, we apply $L^2$ energy method instead of using the estimates of linear problem in the high frequency part. Concerning $L^1$ estimate of solutions we use another type of estimate which has lower singularity at $t=0$ in the high frequency part. Note that by the Korteweg tensor $\rho$ has higher regularity than that of velocity. 

To obtain  Theorem \ref{thm-nonlinear} {\rm (i)} below, the conservation form has an important role similarly to the proof of \cite{Hoff-Zumbrun1} for CNS. 
Note that owing to the smoothing effect of $\rho$ from the Korteweg tensor, even if we consider the conservation form \eqref{CNSK} no derivative loss occurs in the energy method for the high frequency part. The energy estimate is one of key points in the method of \cite{Okita}. Therefore, we have the estimates of solutions to \eqref{CNSK} including Theorem \ref{thm-nonlinear} {\rm (i)} in lower regularity. 
As for CNS, due to the derivative loss of the density, the energy estimate is incompatible with the conservation form. Hence we can not obtain a similar property  to Theorem \ref{thm-nonlinear} {\rm (i)} in lower regularity by using the decomposition method. Indeed, \cite{Hoff-Zumbrun1} shows a similar property to Theorem \ref{thm-nonlinear} {\rm (i)} with an initial data in  $H^{[n/2]+3}$ class for CNS.

This paper is organized as follows. In section 2 notations and lemmas are described which shall be used in this paper. 
In section 3, main results are stated for the linearized and nonlinear problems of \eqref{CNSK} respectively. In section 4, the proof of the time decay estimate of solutions to the linearized problem of \eqref{CNSK} is given. 
In section 4,  
that to nonlinear problem is given.

@

\section{Preliminaries}
In this section we introduce notations which will be used throughout this paper. Furthermore, we introduce some lemmas which will be useful in the proof of the main results.\\

We denote the norm on $X$ by $\|\cdot\|_{X}$ for a given Banach space $X$.

Let $1\leqq p\leqq \infty.$ $L^p$ denotes the usual $L^p$ space on $\mathbb{R}^n$. 
Let $k$ be a nonnegative integer. $W^{k,p}$ and $H^k$ denotes the usual $L^p$ and $L^2$ Sobolev space of order $k$ respectively.  
(As usual, we define that $H^0:=L^2$.) 

For simplicity,  $L^p$ denotes the set of all vector fields 
$w=\trans(w_1, \cdots, w_n)$ on $\mathbb{R}^n$ 
with $w_{j}\in L^p$ $(j=1, \cdots, n)$ 
and $\|\cdot\|_{L^p}$ denotes the norm $\|\cdot\|_{(L^p)^n}$ 
if no confusion will occur. 
Similarly,  a function space $X$ denotes    
the set of all vector fields 
$w=\trans(w_1, \cdots, w_n)$ on $\mathbb{R}^n$ 
with $w_{j}\in L^p$ $(j=1, \cdots, n)$  
and  the norm $\|\cdot\|_{X^n}$ is denoted by $\|\cdot\|_{X}$ if no confusion will occur. 

We take $u=\trans(\phi,m)$ with $\phi\in H^k$ and $m=\trans(m_1,\cdots, m_n)\in H^j$.  
Then the norm $\|u\|_{H^k\times H^j}$ denotes  the norm of $u$ on $H^k\times H^j$, that is, we define
$$
\|u\|_{H^k\times H^j}:=\left(\|\phi\|_{H^k}^2+\|w\|_{H^j}^2\right)^{\frac{1}{2}}.
$$
When $j=k$, for simplicity we denote $H^k\times (H^k)^n$ by $H^k$. 
The norm $\|u\|_{H^k}$ denotes 
the norm $\|u\|_{H^k\times (H^k)^n}$ i.e., we define that 
$$
H^k:=H^k\times (H^k)^n, 
\ \ \ 
\|u\|_{H^k}:=\|u\|_{H^k\times (H^k)^n} 
\ \ \ (u=\trans(\phi,m)).
$$
Similarly, for $u=\trans(\phi,m)\in X\times Y$ with $m=\trans(m_1,\cdots, m_n)$ 
, 
the norm $\|u\|_{X\times Y}$ denotes
$$
\|u\|_{X\times Y}:=\left(\|\phi\|_{X}^2+\|m\|_{Y}^2\right)^{\frac{1}{2}}
\ \ \ (u=\trans(\phi,m)).
$$
If $Y=X^n$, the symbol 
$X$ denotes $X\times X^n$ for simplicity,  
and we define its norm $\|u\|_{X\times X^n}$ by $\|u\|_X$;
$$
X:=X\times X^n, 
\ \ \ 
\|u\|_{X}:=\|u\|_{X\times X^n}
\ \ \ (u=\trans(\phi,m)).
$$


The symbols $\hat{f}$ and $\mathcal{F}[f]$ denote the Fourier transform of $f$ for the space variables; 
\begin{eqnarray}
\hat{f}(\xi)
=\mathcal{F}[f](\xi)
:=\int_{\mathbb{R}^n}f(x)e^{-ix\cdot\xi}dx\quad (\xi\in\mathbb{R}^n).\nonumber
\end{eqnarray}
In addition,  the inverse Fourier transform of $f$ is denoted by $\mathcal{F}^{- 1}[f]$; 
\begin{eqnarray}
\mathcal{F}^{- 1}[f](x)
:=(2\pi)^{-n}\int_{\mathbb{R}^n}f(\xi)e^{i\xi\cdot x}d\xi\quad(x\in\mathbb{R}^n).\nonumber
\end{eqnarray}

The function space $H_{(\infty)}^{k}$ denotes the set of all $u\in H^k$ satisfying $\mbox{supp }\hat{u}\subset\{|\xi|\geq r_{\infty}\}$, where $r_\infty$ denotes a positive constant.  

%
%
%

For operators $L_1$ and $L_2$, we denote by $[L_1,L_2]$ the commutator of $L_1$ and $L_2$, i.e.,  
\begin{eqnarray}
[L_1,L_2]f:=L_1(L_2 f)-L_2(L_1 f).\nonumber
\end{eqnarray}

For a nonnegative number $s$, $[s]$ denotes the integer part of $s$. 

The symbol $``\ast"$ denotes the spatial convolution. 

\vspace{2ex}

We next  state some lemmas which will be used in the proof of the main results.  

\vspace{2ex}

The following lemma is the well-known Sobolev type inequality.

\begin{lem}\label{lem2.1.} Let $s$ satisfy $s > n/2.$ Then there holds the inequality 
\begin{eqnarray}
\|f\|_{L^{\infty}}\leq C\| f\|_{H^{s}}\nonumber
\end{eqnarray}
for $f\in H^{s}.$
\end{lem}


\vspace{2ex}

The following inequalities are stated which are concerned with composite functions.

\vspace{2ex}

\begin{lem}\label{lem2.2.} 
Let $s$ be an integer satisfying $s\geq [n/2]+1$.  
Let $s_{j}$ and $\mu_{(j)}$ ($j=1,\cdots,\ell$) be nonnegative integers 
and multiindices 
satisfying $0\leq |\mu_{(j)}|\leq s_{j}\leq s+|\mu_{(j)}|$, 
$\mu=\mu_{(1)}+\cdots +\mu_{(\ell)}$, 
$s=s_{1}+\cdots+s_{\ell}\geq(\ell-1)s+|\mu|$, respectively. 
Then there holds
\begin{eqnarray}
\parallel \partial^{\mu_{(1)}}_{x}f_{1}\cdots \partial^{\mu_{(\ell)}}_{x}f_{\ell}\parallel_{L^{2}}
\leq C\prod_{1\leq j\leq \ell}\parallel f_{j}\parallel_{H^{s_{j}}}\quad  (f_j \in H^{s_j}).\nonumber
\end{eqnarray}
\end{lem}

See, e.g., \cite{Kagei-Kobayashi} for the proof of Lemma $\ref{lem2.2.}$.

\vspace{2ex}
\begin{lem}\label{lem2.3.} 
Let $s$ be an integer satisfying $s\geq [n/2]+1$. 
Suppose that $F$ is a smooth function on $I$, 
where $I$ is a compact interval of $\mathbb{R}$.  
Then for a multi-index $\alpha$ with $1\leq |\alpha|\leq s$, 
there hold the estimates 
$$
\|[\partial^{\alpha}_{x},F(f_{1})]f_{2}\|_{L^2}
\leq 
C\|F\|_{C^{|\alpha|}(I)}
\left\{1+\|\nabla f_{1}\|_{s-1}^{|\alpha|-1}\right\}
\|\nabla f_{1}\|_{H^{s-1}}\|f_2 \|_{H^{|\alpha|}},
$$ 
for $f_1\in H^s$ with $f_1(x)\in I$ for all $x\in \mathbb{R}^n$ 
and $f_2\in H^{|\alpha|}$; and 
$$
\|[\partial^{\alpha}_{x},F(f_{1})]f_{2}\|_{L^2}
\leq 
C\|F\|_{C^{|\alpha|}(I)}
\left\{1+\|\nabla f_{1}\|_{s-1}^{|\alpha|-1}\right\}
\|\nabla f_{1}\|_{H^{s}}\|f_2 \|_{H^{|\alpha|-1}},
$$
for $f_1\in H^{s+1}$ with $f_1(x)\in I$ for all $x\in \mathbb{R}^n$ 
and $f_2\in H^{|\alpha|-1}$.  
\end{lem}

See, e.g., \cite{Kagei-KawashimaCMP} for the proof of Lemma $\ref{lem2.3.}$.



\section{Main results}
In this section, main results are stated for $(\ref{CNSK})$. \eqref{CNSK} is reformulated as follows.  
Hereafter we assume that $\rho_*=1$ without loss of generality. We set 
$\phi=\rho-1$ and $m=\frac{M}{\gamma}$ where $\gamma=\sqrt{P'(1)}$. Substituting $\phi$ and $m$ into \eqref{CNSK}, then 
we obtain 
\begin{eqnarray}
\left\{
\begin{array}{lll}
\partial_{t}\phi +\gamma\div  m=0,\\
\partial_{t}m-\nu\Delta m-\tilde{\nu}\nabla\div m+\gamma \nabla \phi-\kappa_0 \nabla \Delta\phi=f(u),\\
\phi|_{t=0}=\phi_0, \ \ m|_{t=0}=m_0,
\label{cnsk-nolinear}
\end{array}
\right.
\end{eqnarray}
where $u=\trans(\phi,m)$, $\nu={\mu}$, $\tilde{\nu}={\mu+\mu'}$, $\kappa_0=\frac{\kappa}{\gamma}$, $\phi_0={\rho_0-1}$, $m_0=\frac{M_0}{\gamma}$, 
\begin{eqnarray*}
f(u)&=& -\Big\{\gamma\div(m \otimes m)+\gamma\div(P_{(1)}(\phi)\phi m \otimes m)+\frac{1}{\gamma}\nabla(P_{(2)}(\phi)\phi^2)\\
&&\quad -\nu\Delta(P_{(1)}(\phi)\phi m)-\tilde{\nu}\nabla\div(P_{(1)}(\phi)\phi m)-\div \Phi(\phi)\Big\},\\
P_{(1)}(\phi)
&=&
\int_{0}^{1}f'(1+\tau\phi)d\tau, \ \ f(\tau)=\frac{1}{\tau} \ \ (\tau \in \mathbb{R}),\\
P_{(2)}(\phi)
&=&
\int_{0}^{1}(1-\tau)P''\Big(1+\tau\phi\Big)d\tau,\\
\Phi(\phi)
&=&\kappa_0\Big\{\phi\Delta \phi I_n +(\nabla \phi)\cdot(\nabla \phi)I_n -\frac{|\nabla \phi |^2}{2}I_n-\nabla \phi \otimes \nabla \phi \Big\}.
\end{eqnarray*}

We first consider the time decay estimate of solutions to linearized problem for $(\ref{CNSK})$.  $(\ref{CNSK})$ is  linearized as follows. 
\begin{eqnarray}
\left\{
\begin{array}{lll}
\partial_{t}\phi +\gamma\div  m=0,\\
\partial_{t}m-\nu\Delta m-\tilde{\nu}\nabla\div m+\gamma \nabla \phi-\kappa_0 \nabla \Delta\phi=0,\\
\phi|_{t=0}=\phi_0, \ \ m|_{t=0}=m_0.
\label{cnsk-linear}
\end{array}
\right.
\end{eqnarray}
By taking the Fourier transform of \eqref{cnsk-linear} with respect to the space variable $x$, we obtain the following ordinary differential equation with a parameter $\xi$. 
\begin{eqnarray}
\left\{
\begin{array}{lll}
\partial_{t}\hat{\phi}(t,\xi) +i\gamma\xi\cdot\hat{m}(t,\xi)=0,\\
\partial_{t}\hat{m}(t,\xi)+\nu|\xi|^2 \hat{m}(t,\xi)+\tilde{\nu}\xi(\xi\cdot \hat{m}(t,\xi))+i\gamma \xi \hat{\phi}(t,\xi)
+i\xi \kappa_0 |\xi|^2\hat{\phi}(t,\xi)=0,\\
\hat{\phi}(0,\xi)=\hat{\phi}_0, \ \ \hat{m}(0,\xi)=\hat{m}_0.
\label{cnsk-f}
\end{array}
\right.
\end{eqnarray}
Therefore, the solutions of \eqref{cnsk-linear} are given by the following formulas. Hereafter we define that 
$A=\displaystyle\frac{\nu+\tilde{\nu}}{2}$, $B=\displaystyle\frac{2\gamma}{\nu+\tilde{\nu}}$, $K=\displaystyle\frac{2\sqrt{\kappa_0\gamma}}{\nu+\tilde{\nu}}$. 
If $|\xi|\neq 0,  B/\sqrt{1-K^2}$ when $0< K < 1$ and $|\xi|\neq 0$ when $K \geq 1$, then the Fourier transforms of $\phi$ and $m$ are given by 
\begin{eqnarray}
\hat{\phi}&=&\displaystyle\frac{\lambda_+(\xi) e^{\lambda_-(\xi) t}-\lambda_-(\xi) e^{\lambda_+(\xi) t}}
{\lambda_+(\xi) -\lambda_-(\xi)}\hat{\phi}_0-i\gamma \displaystyle\frac{ e^{\lambda_+(\xi) t}-e^{\lambda_-(\xi) t}}
{\lambda_+(\xi) -\lambda_-(\xi)}\xi\cdot \hat{m}_0,\nonumber\\
\hat{m}&=&e^{-\nu|\xi|^2 t}\hat{m}_0-i\xi(\gamma+\kappa_0|\xi|^2)\left(\displaystyle\frac{ e^{\lambda_+(\xi) t}-e^{\lambda_-(\xi) t}}
{\lambda_+(\xi) -\lambda_-(\xi)}\right)\hat{\phi}_0\nonumber\\
&&\quad +\left(\displaystyle\frac{\lambda_+(\xi) e^{\lambda_+(\xi) t}-\lambda_-(\xi) e^{\lambda_-(\xi) t}}{\lambda_+(\xi) -\lambda_-(\xi)}-e^{-\nu|\xi|^2t}\right)\frac{\xi(\xi\cdot\hat{m}_0)}{|\xi|^2},\label{solution-linear}
\end{eqnarray}
where 
$$
\lambda_{\pm}(\xi)=-A(|\xi|^2\pm\sqrt{|\xi|^4-B^2|\xi|^2-K^2|\xi|^4})
$$
denote roots of the characteristic equation of \eqref{cnsk-f}. Note that when $|\xi|\neq 0$ and $K \geq 1$ it holds that $\lambda_+(\xi) -\lambda_-(\xi) \neq 0$. In addition, due to the term in \eqref{cnsk-linear} from the Korteweg tensor, in contrast to \cite{Kobayashi-Shibata}, the higher order term with respect to $\xi$, i.e.,   
\begin{eqnarray}
i\kappa_0\xi|\xi|^2\left(\displaystyle\frac{ e^{\lambda_+(\xi) t}-e^{\lambda_- (\xi)t}}
{\lambda_+(\xi) -\lambda_-(\xi)}\right)\hat{\phi}_0\label{cnskchara}
\end{eqnarray}
appears in \eqref{solution-linear}.  
On the otherhand, if $0<K< 1$ and $\min\Big\{\frac{1}{2}, \, \displaystyle\frac{B}{2\sqrt{1-K^2}}\Big\}\leq |\xi| \leq 2\displaystyle\frac{B}{\sqrt{1-K^2}}$, then 
$\hat{\phi}$ and  $\hat{m}$ are given by 
\begin{eqnarray}
\hat{\phi}&=&\displaystyle\frac{1}{2\pi i}\oint_{\Gamma}\frac{(z+A|\xi|^2)e^{zt}}
{z^2+(\nu+\tilde{\nu})|\xi|^2 z+\kappa_0\gamma|\xi|^4+\gamma^2|\xi|^2}dz\hat{\phi}_0\nonumber\\
&&\quad-\displaystyle\frac{\gamma}{2\pi }\oint_{\Gamma}\frac{e^{zt}}
{z^2+(\nu+\tilde{\nu})|\xi|^2 z+\kappa_0\gamma|\xi|^4+\gamma^2|\xi|^2}dz\xi\cdot\hat{m}_0,\label{solmidlerho}\\
\hat{m}&=&e^{-\nu|\xi|^2 t}\hat{m}_0-\displaystyle\frac{\gamma \xi}{2\pi }
\oint_{\Gamma}\frac{e^{zt}}
{z^2+(\nu+\tilde{\nu})|\xi|^2 z+\kappa_0\gamma|\xi|^4+\gamma^2|\xi|^2}dz\hat{\phi}_0\nonumber\\
&&\quad+\left(\displaystyle\frac{1}{2\pi i}\oint_{\Gamma}\frac{ze^{zt}}
{z^2+(\nu+\tilde{\nu})|\xi|^2 z+\kappa_0\gamma|\xi|^4+\gamma^2|\xi|^2}dz\right)\frac{\xi(\xi\cdot\hat{m}_0)}{|\xi|^2},\label{solmidlem}
\end{eqnarray}
where $\Gamma$ stands for a closed pass including $\lambda_{\pm}$ and included in the set $\{z\in \mathbb{C}|{\rm Re}z\leq -c_0\}$ and $c_0$ stands for a positive number satisfying that 
$$
\max_{\min\Big\{\frac{1}{2}, \,\frac{B}{2\sqrt{1-K^2}}\Big\}\leq|\xi|\leq 2\frac{B}{\sqrt{1-K^2}}}{\rm Re}\lambda_{\pm}\leq -2c_0.
$$
\vspace{2ex}
We define a cut-off function $\varphi_1$ in $C^{\infty}(\mathbb{R}^n)$ as follows. We assume that $K \neq 1$. 
\begin{eqnarray*}
&&\mbox{When} \ \frac{B}{2\sqrt{|1-K^2|}}> 1,   \ \  \varphi_1(\xi)=
\left\{
\begin{array}{ll}
1 \ \ \mbox{for} \ \   |\xi|\leq \frac{1}{2},\\
0 \ \ \mbox{for} \ \  |\xi|\geq 1.
\end{array}
\right. \\ 
&&\mbox{When} \ \frac{B}{2\sqrt{|1-K^2|}}\leq 1 < \frac{B}{\sqrt{2|1-K^2|}},   \ \ \varphi_1(\xi)=
\left\{
\begin{array}{ll}
1 \ \ \mbox{for} \ \   |\xi|\leq \displaystyle\frac{B}{2\sqrt{|1-K^2|}},\\
0 \ \ \mbox{for} \ \  |\xi|\geq 1.
\end{array}
\right.\\
&&\mbox{When} \ \frac{B}{\sqrt{2|1-K^2|}}\leq 1,  \ \ 
\varphi_1(\xi)=
\left\{
\begin{array}{ll}
1 \ \ \mbox{for} \ \   |\xi|\leq \displaystyle\frac{B}{2\sqrt{|1-K^2|}},\\
0 \ \ \mbox{for} \ \  |\xi|\geq \displaystyle\frac{B}{\sqrt{2|1-K^2|}}.
\end{array}
\right.
\end{eqnarray*}
Furthermore, we define cut-off functions $\varphi_\infty$ and $\varphi_M$ in $C^{\infty}(\mathbb{R}^n)$ by 
\begin{eqnarray*}
&&
\varphi_\infty(\xi)=
\left\{
\begin{array}{ll}
1 \ \  \mbox{for} \ \ |\xi|\geq \displaystyle\frac{2B}{\sqrt{|1-K^2|}},\\
0 \ \  \mbox{for} \ \ |\xi|\leq \displaystyle\frac{\sqrt{2}B}{\sqrt{|1-K^2|}},
\end{array}
\right.\\
&&\varphi_M(\xi)=1-\varphi_1(\xi)-\varphi_\infty(\xi).
\end{eqnarray*}
If $K=1$, $\varphi_1$ and $\varphi_\infty$ are defined as follows. 
\begin{eqnarray*}
&&\varphi_1(\xi)=
\left\{
\begin{array}{ll}
1 \ \ \mbox{for} \ \   |\xi|\leq \frac{1}{2},\\
0 \ \ \mbox{for} \ \  |\xi|\geq 1, 
\end{array}
\right.  \ \ 
\varphi_\infty(\xi)=
\left\{
\begin{array}{ll}
0 \ \ \mbox{for} \ \   |\xi|\leq \frac{1}{2},\\
1 \ \ \mbox{for} \ \  |\xi|\geq 1, 
\end{array}
\right.  \\
&&\varphi_1(\xi)+\varphi_\infty(\xi)=1. 
\end{eqnarray*}
We define solution operators on low frequency part $E_1$ and that on high frequency part $E_\infty$ of \eqref{cnsk-linear} as follows.   
\begin{eqnarray}
&&E_1(t)=(E_{1,\phi}(t),E_{1,m}(t)),\nonumber\\
&&E_\infty(t)=(E_{\infty,\phi}(t),E_{\infty,m}(t)),\label{sol-ope-high}\\
&&E_{1,\phi}(t)(\phi_0,m_0)(x)=\mathcal{F}^{-1}[\varphi_1(\xi)\hat{\phi}(t,\xi)](x),\label{sol-ope-low1}\\
&&E_{1,m}(t)(\phi_0,m_0)(x)=\mathcal{F}^{-1}[\varphi_1(\xi)\hat{m}(t,\xi)](x),\label{sol-ope-low2}\\
&&E_{\infty,\phi}(t)(\phi_0,m_0)(x)=\mathcal{F}^{-1}[(\varphi_M(\xi)+\varphi_\infty(\xi))\hat{\phi}(t,\xi)](x),\label{sol-ope-high1}\\
&&E_{\infty,m}(t)(\phi_0,m_0)(x)=\mathcal{F}^{-1}[(\varphi_M(\xi)+\varphi_\infty(\xi))\hat{m}(t,\xi)](x).\label{sol-ope-high2}
\end{eqnarray}
Note that as in \cite{Tsuda-CNSK} $E_1$ is continuous for $t \geq 0$ on $L^p$ $(1\leq p \leq \infty)$ and $E_\infty$ is continuous for $t \geq 0$ on $H^{s+1}\times H^s$ $(s\geq [n/2]+1)$.   
If an initial time is $t_0>0$ but not $0$ in \eqref{cnsk-linear}, then we write the solution operators by $E_j(t,t_0)$ $(j=1,\infty)$ respectively. 

Concerning the solution $(\phi,m)$ to \eqref{cnsk-linear} and the solution operators on the low frequency part, \eqref{sol-ope-low1} and \eqref{sol-ope-low2}, we obtain the following estimates.

\vspace{2ex}

\begin{thm}\label{est-linear-low} 
{\rm (i)}\quad It holds that for the solution $(\phi,m)$ to \eqref{cnsk-linear} and $t>0$ 
\begin{eqnarray*}
\Big\|\del_t^k\del_x^\alpha \Big\{
\begin{pmatrix}
\phi\\
m
\end{pmatrix}
(t)
-\begin{pmatrix}
0\\
{\cal K}_{\nu}\ast m_{0,in}
\end{pmatrix}
\Big\}
\Big\|_{L^\infty}
\leq C_{k,\alpha,n}(1+t)^{-\left(\frac{3n-1}{4}+\frac{k+|\alpha|}{2}\right)}
[\|\phi_0\|_{L^1}+\|m_0\|_{L^1}],
\end{eqnarray*}
where ${\cal K}_{\nu}={\cal K}_{\nu}(t,x)$ denotes the standard heat kernel and $m_{0,in}$ denotes  a divergence-free part of $m_0$ that are  respectively given by
$$
{\cal K}_{\nu}=\mathcal{F}^{-1}(e^{-\nu|\xi|^2 t}), \ \ m_{0,in}=\mathcal{F}^{-1}\Big\{\Big(I_n-\frac{\xi\trans{\xi}}{|\xi|^2}\Big)m_0\Big\}.
$$

\vspace{1ex}

{\rm (ii)} $($ $L^\infty$ estimate for $E_1$ $)$ For $t>0$ it holds that
\begin{eqnarray*}
\lefteqn{\|\del_t^k\del_x^\alpha E_{1,\phi}(t)(\phi_0,m_0)\|_{L^\infty}}\\
&\leq&C_{k,\alpha,n}(1+t)^{-\left(\frac{3n-1}{4}+\frac{k+|\alpha|}{2}\right)}
[\|\phi_0\|_{L^1}+\|m_0\|_{L^1}];\\[1ex]
\lefteqn{\|\del_t^k\del_x^\alpha E_{1,m}(t)(\phi_0,m_0)\|_{L^\infty}}\\
&\leq&C_{k,\alpha,n}(1+t)^{-\left(\frac{3n-1}{4}+\frac{k+|\alpha|}{2}\right)}
[\|\phi_0\|_{L^1}+\|m_0\|_{L^1}]\\
&&\quad +C_{k,\alpha,n}(1+t)^{-\left(\frac{n}{2}+\frac{k+|\alpha|}{2}\right)}
\|m_0\|_{L^1}.
\end{eqnarray*}

\vspace{1ex}

{\rm (iii)} $($ $L^1$ estimate for $E_1$ $)$ When the space dimension $n\geq 3$ and $n$ is an odd number, then for any $t>0$ we have  the following estimate. 

\begin{eqnarray*}
\lefteqn{\|\del_t^k\del_x^\alpha E_{1}(t)(\phi_0,m_0)\|_{L^1}}\\
&\leq&C_{k,\alpha,n}(1+t)^{\frac{n-1}{4}-\frac{k+|\alpha|}{2}}
[\|\phi_0\|_{L^1}+\|m_0\|_{L^1}]. 
\end{eqnarray*}

\vspace{1ex}

{\rm (iv)} $($ $L^p$-$L^q$ estimate for $E_1$ $)$ For $t>0$ and $1\leq q \leq 2 \leq p \leq \infty$ we have the following estimate.
\begin{eqnarray*}
\lefteqn{\|\del_t^k\del_x^\alpha E_{1}(t)(\phi_0,m_0)\|_{L^p}}\\
&\leq&C_{k,\alpha,n}(1+t)^{-\left(\frac{n}{2}\left(\frac{1}{q}-\frac{1}{p}\right)+\frac{k+|\alpha|}{2}\right)}
[\|\phi_0\|_{L^q}+\|m_0\|_{L^q}].
\end{eqnarray*}
\end{thm}

\vspace{2ex}

\begin{rem}{\rm
Concerning {\rm (iii)}, so far we do not obtain the similar estimate when the space dimension $n\geq 2$ and $n$ is an even number. The key point to obtain {\rm (iii)} is pointwise estimate of the Green function as mentioned in the proof of  Theorem \ref{est-linear-low} below. In the pointwise estimate, we need a great deal of cancellation to overcome the similar difficulty related to the Riesz kernel to that of \cite{Kobayashi-Shibata}. As in \cite{Hoff-Zumbrun2} it seems that the pointwise estimates of the Green function  are different between odd dimensional case and even dimensional case by the Huygens principle and the estimate with even dimensional case is more complicated than that of odd dimensional case. Hence, more delicate analysis is needed to obtain the similar estimate to {\rm (iii)} in even dimensional case.  Note that diffusion wave property, especially the retardation of the parabolic decay in $L^p$ ($p <2$) occurs for multi dimensional case as in  \cite{Hoff-Zumbrun1}.  
}\end{rem}

\vspace{2ex}

\begin{rem}
{\rm We discuss the optimality of the decay exponents in Theorem \ref{est-linear-low} below. Concerning {\rm (i)} and {\rm (iv)}, the first approximation of solutions is ${\cal K}_{\nu}\ast m_{0,in}$, that is, the Stokes flow part and optimality of decay exponents of solutions to the Stokes equation is well known. Concerning {\rm (ii)} and {\rm (iii)}, Hoff and Zumbrun \cite{Hoff-Zumbrun2} consider some linear artificial viscous equation whose solutions approximate behavior of solutions to the linearized compressible Navier-Stokes equation. They give not only upper bounds but also lower ones of the Green function and verify that  decay exponents of $L^p$ bounds are sharp at least up to a logarithmic term. The decay exponents coincide with those of {\rm (ii)} and {\rm (iii)}. Therefore we think that the decay rates of {\rm (ii)} and {\rm (iii)} are optimal. }
\end{rem}

\vspace{2ex}

Furthermore, the following estimate holds for the solution operator on the high frequency part \eqref{sol-ope-high}.

\vspace{2ex}

\begin{thm}\label{est-linear-high}  $($ $L^p$-$L^p$ estimate for $E_\infty$ $)$
Let $1\leq p \leq \infty$. Then it holds that
\begin{eqnarray}
\lefteqn{\|\del_t^k \del_x^\alpha E_{\infty}(t)(\phi_0,m_0)\|_{L^p}}\nonumber\\
&\leq & C_{k,\alpha, n}e^{-ct} \Big\{(1+t^{-\delta_1-\frac{|\alpha|}{2}-k})\Big[\|\phi_0\|_{L^p}+\|m_0\|_{L^p}\Big]+
t^{-\delta_2 -\frac{|\alpha|}{2}-k}\|\phi_0\|_{L^p}\Big\}\label{est-linear-high1}
\end{eqnarray}
for $t>0$, $k \geq 0$ and $|\alpha|\geq 0$,  
where $(\delta_1,\delta_2 )=(1/2,1),(1,3/2)$ for $K\neq 1$ and $K=1$ respectively. In addition, when $K=1$ and $0<t \leq1$, we have the following estimate. 
\begin{eqnarray}
\|E_{\infty}(t)(\phi_0,m_0)\|_{L^p}\leq C_{n}t^{-\frac{1}{2}-\sigma_0}\Big[\|\phi_0\|_{L^p}+\|m_0\|_{L^p}\Big]+t^{-\delta_2 }\|\phi_0\|_{L^p}\label{est-linear-high2}
\end{eqnarray}
for $1\leq p \leq \infty$, where $\sigma_0$ is any positive number satisfying $0< \sigma_0 <1/2$. 
\end{thm}

\vspace{2ex}

\begin{rem}{\rm
Theorem \ref{est-linear-high} implies smoothing effect of solutions to the linearized problem in the high frequency part. The estimate \eqref{est-linear-high2} has lower singularity at $t=0$ than \eqref{est-linear-high1} with $K=1$. }
\end{rem}

\vspace{2ex}

We next consider time decay estimates of solutions to the nonlinear problem \eqref{CNSK}. 
The following $L^\infty$, $L^2$ and $L^1$ estimates are stated for the solution $u=\trans(\phi,m)$ to the system \eqref{cnsk-nolinear}. 

\vspace{2ex}

\begin{thm}\label{thm-nonlinear} Let $u=\trans(\phi,m)$ be the solution to \eqref{cnsk-nolinear}. 

\noindent{\rm (i)}  Let $E=E(t)$ be the solution operator for the linearized problem \eqref{cnsk-linear} defined by $E=E_1+E_{\infty}$. We assume that $u_0=\trans(\phi_0,m_0)\in (H^{s+1}\times H^s) \cap L^1$, where $s$ denotes a nonnegative integer satisfying $s\geq [n/2]+1$. We define the norm $|||u_0|||_s$ by 
$$
|||u_0|||_s=\|u_0\|_{(H^{s+1}\times H^s ) \cap L^1}. 
$$
There exists a constant $\epsilon_1>0$ such that if 
$|||u_0|||_s\leq \epsilon_1$, 
then for $t\geq 0$ we have the estimate 
\begin{eqnarray*}
 \left\|
\trans(\phi,m)
(t)
- E(t)
(\phi_0,m_0)
\right\|_{L^\infty}
 \leq C(1+t)^{-\frac{n}{2}-\frac{1}{2}}\delta_1(t) |||u_0|||_s,
\end{eqnarray*} 
where $\delta_1(t)=1$ for $n\geq 3$ and $\delta_1(t)=\log (1+t)$ for $n=2$. 
\vspace{2ex}

{\rm (ii)} Under the assumption of {\rm (i)} it holds that 
\begin{eqnarray*}
&&\|\phi(t)\|_{L^\infty}\leq C(1+t)^{q(n)}|||u_0|||_s, \\
&&\|m(t)\|_{L^\infty}\leq C(1+t)^{-\frac{n}{2}}|||u_0|||_s,
\end{eqnarray*}
where $q(n)=-\frac{3n-1}{4}$ for $n=2,3$ and $q(n)=-\frac{n}{2}-\frac{1}{2}$ for $n\geq 4$. 
Furthermore, it also holds that for $t\geq 0$ and $k=0,1$  
\begin{eqnarray*}
\|\nabla ^k u(t)\|_{L^2}\leq C(1+t)^{-\frac{n}{4}-\frac{k}{2}}|||u_0|||_s.
\end{eqnarray*}
\vspace{2ex}

{\rm (iii)} Let $n$ be an odd number satisfying $n\geq 3$. There exists a constant $\epsilon_2>0$ such that if 
$|||u_0|||_s\leq \epsilon_2$, then the following estimate is true for $t\geq 1$.  
\begin{eqnarray*}
\|u(t)\|_{L^1}\leq Ct^{\frac{n-1}{4}}|||u_0|||_s. 
\end{eqnarray*}
\end{thm}

\vspace{2ex}

\begin{rem}{\rm 
Concerning the first approximation of $u$, i.e., ${\cal K}_{\nu}\ast m_{0,in}$, we have the following estimate. 
$$
\|{\cal K}_{\nu}\ast m_{0,in}\|_{L^\infty} \leq C(1+t)^{-\frac{n}{2}}|||u_0|||_s. 
$$
}
\end{rem}

\vspace{2ex}

\begin{rem}{\rm
In Theorem \ref{thm-nonlinear}, the diffusion wave property of the solution appears in $L^\infty$  and $L^1$ estimates.} 
\end{rem}

\vspace{2ex}

\section{Proof of the estimates for solution to the linear problem} 

In this section, we give the proofs of Theorem \ref{est-linear-low} and Theorem \ref{est-linear-high}. To prove Theorem \ref{est-linear-low}, we put 
\begin{eqnarray*}
L_{11,j}(t,x)&=&\mathcal{F}^{- 1}\Big\{\displaystyle\frac{\lambda_+(\xi) e^{\lambda_-(\xi) t}-\lambda_-(\xi) e^{\lambda_+(\xi) t}}{\lambda_+(\xi) -\lambda_-(\xi)}\varphi_j(\xi)\Big\}(x), \\
L_{12,j}(t,x)&=&\mathcal{F}^{- 1}(-i\gamma\trans\xi\hat{L}_j),\\
\hat{L}_j(t,\xi)&=&\displaystyle\frac{ e^{\lambda_+(\xi) t}-e^{\lambda_-(\xi) t}}{\lambda_+(\xi) -\lambda_-(\xi)}\varphi_j(\xi), \\
L_{21,j}(t,x)&=&\mathcal{F}^{-1}\{-\xi(i\gamma +k_0|\xi|^2 )\hat{L}_j\},\\
L_{22,j}(t,x)&=&K_{1,j}(t,x)+K_{2,j}(t,x)-K_{3,j}(t,x),\\
K_{1,j}(t,x)&=&\mathcal{F}^{- 1}\Big[e^{-\nu|\xi|^2 t}\varphi_j(\xi)\Big](x)I_n,\\
K_{2,j}(t,x)&=&\mathcal{F}^{- 1}\Big\{\displaystyle\frac{\lambda_+(\xi) e^{\lambda_+(\xi) t}-\lambda_-(\xi) e^{\lambda_-(\xi) t}}{\lambda_+(\xi) -\lambda_-(\xi)}\frac{\xi \trans{\xi}}{|\xi|^2}\varphi_j(\xi)\Big\}(x), \\
K_{3,j}(t,x)&=&\mathcal{F}^{- 1}\Big[e^{-\nu|\xi|^2 t}\frac{\xi \trans{\xi}}{|\xi|^2}\varphi_j(\xi)\Big](x)
\end{eqnarray*}
for $j=1,\infty$. We see from $\eqref{solution-linear}$ that 
\begin{eqnarray}
E_j(t)(\phi_0,m_0)=
\begin{pmatrix}
L_{11,j}(t,\cdot) & L_{12,j}(t,\cdot)\\
L_{21,j}(t,\cdot) & L_{22,j}(t,\cdot)\\
\end{pmatrix}
*
\begin{pmatrix}
\phi_0\\
m_0
\end{pmatrix}
\label{convolution-solutionfomula}
\end{eqnarray}
for $j=1, \infty$. We first show Theorem \ref{est-linear-low} {\rm (ii)}. 
Concerning the proof of Theorem \ref{est-linear-low} {\rm (ii)}, it is enough to show the following proposition by the same reason as that in \cite[Theorem 2.1 (1)]{Kobayashi-Shibata} based on the Young inequality. 

\vspace{2ex}

\begin{prop}\label{low-kernelest}
We set  
$$
K_{\psi}=\mathcal{F}^{- 1}\Big[\frac{e^{\lambda_+(\xi) t}-e^{\lambda_-(\xi) t}}{\lambda_+(\xi) -\lambda_-(\xi)}\psi(\xi)\varphi_1(\xi)\Big],
$$
where $\psi=\psi(\omega)\in C^\infty(S^{n-1})$, $S^{n-1}=\{\xi \in \mathbb{R}^n; |\xi|=1\}$ and $\psi(\xi)=\psi(\xi/|\xi|)$. Then it holds that 
$$
\|\del_t^k \del_x^\alpha K_\psi(t,\cdot)\|_{L^\infty}\leq C_{k,\alpha,n}(1+t)^{-\left(\frac{3n-3}{4}+\frac{k+|\alpha|}{2}\right)}
$$
for $n \geq 2$ and $t> 0$. 
\end{prop}

\vspace{2ex}

Proposition \ref{low-kernelest} is yielded as follows; As for  the estimate near a light cone, that is, for $\{(t,x); |x|\geq R_0 t, t\geq \max(1,R/R_0)^4\}$, using the stationary phase method as that in \cite{Kobayashi-Shibata} we directly obtain Proposition \ref{low-kernelest}, where $R$ and $R_0$ are some suitable positive constants used in the argument. 
In  the case such that $t\geq 1$ and $|x|\leq R_0 t$, we set 
$$
D=\Big|\frac{(\nu+\tilde{\nu})^2-4\kappa_0\gamma}{4\gamma^2}\Big|.
$$
If $(\nu+\tilde{\nu})^2 <4\kappa_0\gamma$, we define $f_1$ and $g_1$ by
\begin{eqnarray*}
f_1(\xi)=1+|\xi|^2 g_1(|\xi|^2), \ g_1(s)=\frac{D}{2}\displaystyle\int_0^1\frac{1}{\sqrt{1+\theta s D}}d\theta.
\end{eqnarray*}
If $(\nu+\tilde{\nu})^2 \geq 4\kappa_0\gamma$, we define $f_2$ and $g_2$ by
\begin{eqnarray*}
f_2(\xi)=1+|\xi|^2 g_2(|\xi|^2), \ g_2(s)=-\frac{D}{2}\displaystyle\int_0^1\frac{1}{\sqrt{1+\theta s D}}d\theta.
\end{eqnarray*}
Note that  
\begin{eqnarray}
e^{\lambda_{\pm}(\xi)t}=e^{-A|\xi|^2 t} e^{\mp i \gamma |\xi|f_j(|\xi|)t}  \ \ (j=1,2). \label{diffusion-part}
\end{eqnarray}
In contrast to \cite{Kobayashi-Shibata}, new terms $\kappa_0\nabla \Delta \phi$ and \eqref{cnskchara} appear in the linearized problem and the solution formula respectively. 
However, $f_j$ and $g_j$ $(j=1,2)$ have the same order of $|\xi|$ as $f$ and $g$ used in the proof of \cite[Theorem 2.1 (1)]{Kobayashi-Shibata} as $|\xi|$ goes to $0$. 
Therefore a similar manner to the proof of \cite[Theorem 2.1 (1)]{Kobayashi-Shibata} based on \eqref{diffusion-part}, the well-known formulas for fundamental solution to wave equation and Shimizu and Shibata \cite[Theorem 2.3]{Shibata-Shimizu} shows Proposition \ref{low-kernelest} in  the case such that $t\geq 1$ and $|x|\leq R_0 t$.   
Since direct calculation shows that 
$$
|\del_t^k \del_x^\beta K_\psi(t,x)|\leq C_{k,\beta,n}  \ \ \mbox{for} \ \ 0\leq t \leq \max(1,R/R_0)^4, 
$$ 
we get Proposition \ref{low-kernelest}. 
Using $f_j$ and $g_j$ $(j=1,2)$, Theorem \ref{est-linear-low} {\rm (iii)} and Theorem \ref{est-linear-low} {\rm (iv)} are directly verified by a similar proof to that of \cite[Theorem 2.1 (2), Theorem 2.3]{Kobayashi-Shibata}.

Theorem \ref{est-linear-low} {\rm (i)} easily follows from Proposition \ref{low-kernelest}, definitions of $L_{i,j}$ and the Young inequality.

\vspace{2ex}

We next show Theorem \ref{est-linear-high}. Similarly to \cite{Kobayashi-Shibata}, we  define $K_{-,\infty}(t)$ and $K_{-,\ell}(t,x)$ by 
$$
K_{-,\infty}(t)m_0(x)=\mathcal{F}^{- 1}\Big[\frac{\lambda_- (\xi) e^{\lambda_-(\xi) t}}
{\lambda_+(\xi)-\lambda_-(\xi)}
\frac{\xi\trans{\xi}}{|\xi|^2}\varphi_\infty(\xi)\hat{m}_0(\xi)\Big](x), 
$$
and 
$$
K_{-,\ell}(t,x)=\mathcal{F}^{- 1}\Big[\frac{\lambda_-^\ell (\xi)e^{\lambda_-(\xi) t}}{\lambda_+(\xi)-\lambda_-(\xi)} 
\frac{\xi\trans{\xi}}{|\xi|^2}\varphi_\infty(\xi)\Big](x)
$$
for $\ell\geq  0$. Note that 
$$
\del_{t}^\ell\del_x^\alpha K_{-,\infty}(t) m_0 =\del_x^\alpha K_{-,\ell+1}(t,\cdot) *m_0.
$$
We put $L_{\pm}(t)$, $M_{\pm, \beta}(t)$, $K_{+,\infty}(t)$ and $K_{1,\infty}(t)$ by the same forms as those used in  \cite{Kobayashi-Shibata}, 
i.e., we put 
\begin{eqnarray*}
L_{\pm}(t)u_0(x)&=&\mathcal{F}^{- 1}\Big[\frac{\lambda_{\mp} (\xi) e^{\lambda_{\pm}(\xi) t}}
{\lambda_+(\xi)-\lambda_-(\xi)}
\varphi_\infty(\xi)\hat{u}_0(\xi)\Big](x),\\
M_{\pm, \beta}(t)u_0(x)&=&\mathcal{F}^{- 1}\Big[\frac{\xi^{\beta} e^{\lambda_{\pm}(\xi) t}}
{\lambda_+(\xi)-\lambda_-(\xi)}
\varphi_\infty(\xi)\hat{u}_0(\xi)\Big](x),  \ \ |\beta|=1, \\
K_{+,\infty}(t)m_0(x)&=&\mathcal{F}^{- 1}\Big[\frac{\lambda_+ (\xi) e^{\lambda_+(\xi) t}}
{\lambda_+(\xi)-\lambda_-(\xi)}
\frac{\xi\trans{\xi}}{|\xi|^2}\varphi_\infty(\xi)\hat{m}_0(\xi)\Big](x),\\
K_{1,\infty}(t)m_0(x)&=&\mathcal{F}^{- 1}\Big[e^{-\nu|\xi|^2 t}\Big(I_n-
\frac{\xi\trans{\xi}}{|\xi|^2}\Big)\varphi_\infty(\xi)\hat{m}_0(\xi)\Big](x). 
\end{eqnarray*}
Note that $L_{ab,\infty}$ $(a=1,2, b=1,2)$ and $K_{a,\infty}$ $(a=1,2,3)$ are linear combinations by $L_{\pm}(t)$, $M_{\pm, \beta}(t)$, $\Delta M_{\pm, \beta}(t)$, $K_{\pm,\infty}(t)$ and $K_{1,\infty}(t)$.  
We shall show that 
\begin{eqnarray}
\|\del_{t}^\ell \del_x^{\alpha}K_{-,\infty}(t) m_0\|_{L^p}\leq Ct^{-\frac{1}{2}-\frac{|\alpha|}{2}-\ell}e^{-Ct}
\|m_0\|_{L^p}\label{est Kinf}
\end{eqnarray}
for $K\neq 1$, $t> 0$, $\ell \geq 0$, $|\alpha|\geq 0$ and $1\leq p \leq \infty$. When $K=1$, we also show that 
\begin{eqnarray}
\|\del_{t}^\ell \del_x^{\alpha}K_{-,\infty}(t) m_0\|_{L^p}\leq Ct^{-1-\frac{|\alpha|}{2}-\ell}e^{-Ct}
\|m_0\|_{L^p}.\label{est Kinf2}
\end{eqnarray}
Indeed, if $K \neq 1$, since $\kappa_0>0$, we have that 
\begin{eqnarray}
\lambda_- =
\left\{
\begin{array}{ll}
-A|\xi|^2+A|1-K^2|^{\frac{1}{2}}|\xi|^2\sqrt{1-\frac{B^2}{|1-K^2||\xi|^2}} \ \ (\nu+\tilde{\nu}> 2\sqrt{\kappa_0 \gamma}),\\
-A|\xi|^2+i A|1 -K^2|^{\frac{1}{2}}|\xi|^2\sqrt{1-\frac{B^2}{|1-K^2||\xi|^2}} \ \ (\nu+\tilde{\nu} < 2\sqrt{\kappa_0 \gamma}).\label{lamda-}
\end{array}
\right.
\end{eqnarray}
If $K =1$, we have that 
\begin{eqnarray}
\lambda_- = -A|\xi|^2 + i\gamma |\xi|. \label{lamda-2}
\end{eqnarray}
Hence, if $\nu+\tilde{\nu}> 2\sqrt{\kappa_0 \gamma}$, we can rewrite $\lambda_-$ to 
\begin{eqnarray*}
\lambda_-=-A(1-|1-K^2|^{\frac{1}{2}})|\xi|^2-\frac{A B^2}{2|1-K^2|^{\frac{1}{2}}}
-\frac{AB^4}{8|1-K^2|^{\frac{3}{2}}|\xi|^2}g_3\Big(\frac{B^2}{|1-K^2||\xi|^2}\Big),
\end{eqnarray*}
where 
$$
g_3 (s)=\displaystyle\int_{0}^{1}(1-\theta s)^{-\frac{3}{2}}d\theta, 
$$
 and note that 
$$
 1-|1-K^2|^{\frac{1}{2}}=1-\Big|1-\frac{4 \kappa_0\gamma}{(\nu+\tilde{\nu})^2}\Big|^{\frac{1}{2}}> 0.
$$
Furthermore, if $\nu+\tilde{\nu} \leq 2\sqrt{\kappa_0 \gamma}$ we see that ${\rm Re}\lambda_-=-A|\xi|^2$. 
Hence we see from \eqref{lamda-} and \eqref{lamda-2}  that for $|\xi|\geq \frac{2B}{\sqrt{|1-K^2|}}$ with  $K \neq 1$ or $|\xi|\geq 1$ with $K =1$ and $|\beta|\geq 0$ there exist positive constants $c_1$ and $c_2$ such that 
\begin{eqnarray*}
\left|\del_{\xi}^\beta \left\{\lambda_-^{\ell+1} (\xi) e^{\lambda_-(\xi) t}
\frac{\xi \trans{\xi}}{|\xi|^2}\varphi_\infty(\xi)\right\}\right| 
\leq C_{\ell,\beta}t^{-\ell}e^{-c_1t} |\xi|^{2-|\beta|}e^{-c_2|\xi|^2 t}
\end{eqnarray*}
In addition, we see that 
\begin{eqnarray}
\frac{1}{|\lambda_+(\xi)-\lambda_-(\xi)|} \leq 
\left\{
\begin{array}{ll}
C|\xi|^{-2} \ \ (K \neq 1,  \ |\xi|\geq \frac{2B}{\sqrt{|1-K^2|}}),\\
C|\xi|^{-1} \ \ (K = 1,  \ |\xi|\geq 1).
\end{array}
\right.
\label{lamdahikizan} 
\end{eqnarray}
Therefore, using
$$
e^{ix\cdot \xi}= \sum_{j=1}^{n}\displaystyle\frac{x_j}{i|x|^2}\del_{\xi_j}e^{ix\cdot \xi}
$$
and the integration by parts $n+1$ times, we obtain that for $K \neq 1$ 
\begin{eqnarray}
|\del_x^{\alpha}K_{-,\ell+1}(t,x)|&\leq& 
C\displaystyle\frac{t^{-\ell}e^{-c_1 t}}{|x|^{n+1}}
\int_{|\xi|\geq \frac{2B}{\sqrt{|1-K^2|}} } |\xi|^{-n-1+|\alpha|}e^{-c_2|\xi|^2 t }
d\xi\nonumber\\
&& \leq C\displaystyle\frac{t^{-\ell-\frac{|\alpha|}{2}}e^{-c_1 t}}{|x|^{n+1}}.\label{estsmallk inf}
\end{eqnarray}
On the other hand, we also obtain from \eqref{lamda-} that 
\begin{eqnarray}
|\del_x^{\alpha}K_{-,\ell+1}(t,x)|&\leq& 
C\displaystyle t^{-\ell}e^{-c_1 t}
\int_{|\xi|\geq \frac{2B}{\sqrt{|1-K^2|}} } |\xi|^{|\alpha|}e^{-c_2|\xi|^2 t }
d\xi\nonumber\\
&& \leq C\displaystyle t^{-\ell-\frac{|\alpha|+n+1}{2}}e^{-c_1 t}.\label{estsmallk inf-2}
\end{eqnarray}
It follows from  \eqref{estsmallk inf} and \eqref{estsmallk inf-2} that 
 \begin{eqnarray}
\|\del_x^{\alpha}K_{-,\ell}(t,\cdot)\|_{L^1}&\leq& C_{\ell, \alpha,n}t^{-\ell-\frac{|\alpha|}{2}}e^{-c_1 t}
\left\{t^{-\frac{n+1}{2}}\displaystyle\int_{|x|\leq \sqrt{t}}dx+\displaystyle\int_{|x|\geq \sqrt{t}}\frac{1}{|x|^{n+1}}dx\right\}\nonumber\\
&\leq & C_{\ell,\alpha,n}t^{-\frac{1}{2}-\frac{|\alpha|}{2}-\ell}e^{-c_1t}. \label{kmainasuest}
\end{eqnarray}
\eqref{kmainasuest} and the Young inequality derive \eqref{est Kinf}. When $K=1$, \eqref{est Kinf2} is verified by a similar argument to the proof of \eqref{est Kinf}. 
Furthermore, when $K=1$, we see that for $0< t \leq1$ and $0<\sigma_0<1/2$ 
\begin{eqnarray}
\Big|\del_\xi^{\alpha} \mathcal{F}(K_{-,1})(t,\xi)\Big| \leq C|\xi|^{-|\alpha|-2\sigma_0}t^{-\frac{1}{2}-\sigma_0}.\label{t0hukin} 
\end{eqnarray}
\eqref{t0hukin} and the strong $L^p$ multiplier theorem (\cite[Proposition 4.2]{Hoff-Zumbrun1}) show that 
\begin{eqnarray}
\|K_{-,\infty}(t,x)m_0\|_{L^p} \leq C_{n}t^{-\frac{1}{2}-\sigma_0}\|m_0\|_{L^p}.\label{Kinf0hukin} 
\end{eqnarray}
for $1\leq p \leq \infty$ and $0<t \leq 1$. 
$L_{\pm}(t)$, $M_{\pm, \beta}(t)$, $K_{+,\infty}(t)$ and $K_{1,\infty}(t)$ are estimated similarly to  \eqref{est Kinf},  \eqref{est Kinf2} and \eqref{Kinf0hukin}.  In addition, we see from the estimate of $M_{\pm, \beta}(t)$ and \eqref{lamdahikizan} that 
$$
\|\del_t^{\ell}\del_{x}^{\alpha}\Delta M_{\pm, \beta}(t)\|_{L^1}\leq C t^{-\delta_2-\frac{|\alpha|}{2}-\ell}e^{-c_1t}, 
$$
where $\delta_2$ is the same one as that in Theorem \ref{est-linear-high}.    
Concerning  the estimate in the middle frequency part the desired estimate  directly follows from the solution formulas \eqref{solmidlerho} and \eqref{solmidlem} and the same manner as that in \cite{Kobayashi-Shibata}, i.e., we have the estimate 
\begin{eqnarray}
\|\del_t^k \del_x^{\alpha}{\mathcal F}^{-1}[\varphi_{M}(\xi)(\hat{\phi},\hat{m})(t,\xi)]\|_{L^p}\leq C_{k,\alpha,n}
e^{-ct}\|\trans(\phi_0, m_0)\|_{L^p}\label{midleest-high}
\end{eqnarray}
for $1\leq p \leq \infty$. 
This together with estimates $L_{\pm}(t)$, $M_{\pm, \beta}(t)$, $K_{\pm,\infty}(t)$ and $K_{1,\infty}(t)$ shows 
Theorem \ref{est-linear-high}.  This completes the proof. $\hfill\square$

\section{Proof of the estimates for solution to the nonlinear problem} 

In this section, we give the proof of Theorem \ref{thm-nonlinear}. Concerning \eqref{cnsk-nolinear}, 
set  
\begin{eqnarray*}
u=\trans(\phi, m), \ \ u_0=\trans(\phi_0,m_0), \ \ A=
\begin{pmatrix}
0 & \gamma\div \\
\gamma \nabla & -\nu \Delta -\tilde{\nu}\nabla \div -\kappa_0 \nabla \Delta
\end{pmatrix}
.
\end{eqnarray*}
Then \eqref{cnsk-nolinear} is rewritten as 
\begin{eqnarray}
\del_t u+Au=F(u),  \ \ u|_{t=0}=u_0,\label{nonlinear-U}
\end{eqnarray}
where $F(u)=\trans(0,f(u))$. Note that we already have existence of global $H^{s+1}\times H^s$ solutions in the introduction.  Hence our task is to discuss a priori estimates of the solutions to \eqref{nonlinear-U}. 

We consider decomposition of solutions to \eqref{nonlinear-U} into low frequency part and high frequency part as in \cite{Okita}. Operators $P_1$ and $P_\infty$ on $L^2$ are defined as follows for the decomposition;
$$
P_1 u=\mathcal{F}^{-1}(\varphi_{1} \hat{u}),  \ \ P_\infty u=\mathcal{F}^{-1}((\varphi_{M}+\varphi_{\infty}) \hat{u}). 
$$  
We have the following properties for $P_j$ $(j=1, \infty)$. 

\vspace{2ex}

\begin{lem}\label{lemP_1}{\rm \cite[Lemma 4.2]{Okita}}
Let $k$ be a nonnegative integer. 
Then $P_{1}$ is a bounded linear operator from $L^2$ to $H^{k}$. 
In fact, it holds that 
\begin{eqnarray}
\|\nabla^{k}P_{1}f\|_{L^2}\leq C\|f\|_{L^2}\qquad (f\in L^{2}).\nonumber
\end{eqnarray}
As a result, for any $2\leq p\leq \infty$, $P_1$ is bounded from $L^2$ to $L^p$. 
\end{lem}
\vspace{2ex}

\begin{lem}\label{lemPinfty}{\rm \cite[Lemma 4.2]{Okita}, \cite[Lemma 4.4]{Kagei-Tsuda}}
{\rm (i)} 
Let $k$ be a nonnegative integer. 
Then $P_\infty$ is a bounded linear operator on $H^k$. 

\vspace{1ex}
{\rm (ii)} 
There hold the inequalities 
\begin{eqnarray*}
\|P_\infty f\|_{L^2}
& \leq &
C\|\nabla f\|_{L^2} 
\ \ (f\in H^1),
\\[2ex]
\|F_{(\infty)}\|_{L^2}
& \leq & 
C\|\nabla F_{(\infty)}\|_{L^2} 
\ \ (F_{(\infty)}\in H^{1}_{(\infty)}).
\end{eqnarray*}
\end{lem}

\vspace{2ex}

Let $u=\trans(\phi, m)$ be a solution to \eqref{nonlinear-U}. Based on the Duhamel principle, similarly to \cite[Proposition 4.3]{Okita}, we decompose the solution to \eqref{nonlinear-U} into low frequency part and high frequency part as  
\begin{eqnarray}
u_j(t)=E_j(t) u_{0j}+ \int_{0}^{t}E_j(t,\tau)P_j F(u_1+u_\infty)(\tau) d\tau\label{low}
\end{eqnarray}
where $u_j=P_j u$ and $u_{0j}=P_j u_0$ $(j=1,\infty)$. 

\vspace{2ex}
We prove Theorem \ref{thm-nonlinear} {\rm (ii)} before we prove {\rm (i)}. We set
\begin{eqnarray*}
M_1(t)=\sup_{0\leq \tau \leq t}\{(1+\tau)^{q(n)}\|\phi_1\|_{L^\infty}
+(1+\tau)^{\frac{n}{2}}\|m_1\|_{L^\infty}+
\sum_{j=0}^{1}(1+\tau)^{\frac{n}{4}+\frac{j}{2}}\|\nabla^j u_1\|_{L^2}\},
\end{eqnarray*}
where $\phi_1=P_1 \phi$ and $m_1=P_1m$.  Furthermore, we also set 
\begin{eqnarray*}
M_\infty(t)=\sup_{0\leq \tau \leq t}(1+\tau)^{\frac{3n-1}{4}}\|u_\infty\|_{H^{s+1}\times H^s}  
\end{eqnarray*}  
for the space dimension $n=2,3$ and 
\begin{eqnarray*}
M_\infty(t)=\sup_{0\leq \tau \leq t}(1+\tau)^{\frac{n}{2}+\frac{1}{2}}\|u_\infty\|_{H^{s+1}\times H^s}  
\end{eqnarray*} 
for the space dimension $n\geq 4$.  
By Theorem \ref{est-linear-low} {\rm (ii)} and \eqref{low}, it holds that 
\begin{align}
&\| u_1(t)\|_{L^\infty}\leq \| E_1(t)u_{01}\|_{L^\infty} +\displaystyle
\int_{0}^{t}\| E_1(t,\tau)P_1 F(u_1+u_\infty)\|_{L^\infty}d\tau,\label{est-linear-lowinfty}\\
&\| E_{1,\phi}(t)u_{01}\|_{L^\infty}
\leq  C(1+t)^{-\frac{3n-1}{4}}\|u_{01}\|_{L^1}, \label{est-linear-lowinfty2}\\
&\| E_{1,m}(t)u_{01}\|_{L^\infty}
\leq  C(1+t)^{-\frac{n}{2}}\|u_{01}\|_{L^1}.\label{est-linear-lowinfty3}
\end{align}
Concerning estimate of the second term of right hand side in \eqref{est-linear-lowinfty}, 
we obtain that 
\begin{eqnarray}
\displaystyle
\int_{0}^{t}\|E_{1,\phi}(t,\tau)P_1 F(u_1+u_\infty)\|_{L^\infty}d\tau\leq C(1+t)^{q(n)}M^2 (t),\label{nonlinearityest-low1}
\end{eqnarray}
where we define that $M(t)=M_1(t)+M_\infty (t)$.  
In fact, as for the convection term $P_1\div(m \otimes m)$ in the nonlinearity $P_1 F(u_1+u_\infty)$, we see from Lemma \ref{lemP_1} that 
\begin{eqnarray*}
\|P_1\div(m \otimes m)\|_{L^1}\leq C\|m\|_{L^2}\|\nabla m\|_{L^2}\leq C(1+t)^{-\frac{n}{2}-\frac{1}{2}}M^2(t)
\end{eqnarray*}
and then it follow from Theorem \ref{est-linear-low} {\rm (ii)} that
\begin{eqnarray*}
\displaystyle
\int_{0}^{t}\Big\|E_{1,\phi}(t,\tau)P_1 
\begin{pmatrix}
0 \\
\div(m \otimes m)
\end{pmatrix}
\Big\|_{L^\infty}d\tau
&\leq &C\displaystyle
\int_{0}^{t}(1+t-\tau)^{-\frac{3n-1}{4}}(1+\tau)^{-\frac{n}{2}-\frac{1}{2}}d\tau M^2 (t)\\
&\leq & C(1+t)^{q(n)}M^2 (t).
\end{eqnarray*}
Since other nonlinear terms are estimated similarly, we have \eqref{nonlinearityest-low1}. 
We next show that
\begin{eqnarray}
\displaystyle
\int_{0}^{t}\| E_{1,m}(t,\tau)P_1 F(u_1+u_\infty)
\|_{L^\infty}d\tau
&\leq &C(1+t)^{-\frac{n}{2}-\frac{1}{2}}\delta_1(t) M^2 (t),\label{m-partest}
\end{eqnarray}
where $\delta_1(t)$ is the same one in Theorem \ref{thm-nonlinear} {\rm (i)}. 
 Indeed, 
we define $I_1$ and $I_2$ by
$$
I_1 = \displaystyle\int_{0}^{\frac{t}{2}}\Big\| E_{1,m}(t,\tau)P_1 
\begin{pmatrix}
0 \\
\div(m \otimes m)
\end{pmatrix}
\Big\|_{L^\infty}d\tau
$$
and 
$$
I_2 = \displaystyle\int_{\frac{t}{2}}^{t}\Big\| E_{1,m}(t,\tau)P_1 
\begin{pmatrix}
0 \\
\div(m \otimes m)
\end{pmatrix}
\Big\|_{L^\infty}d\tau.
$$
$I_1$ is estimated by Theorem \ref{est-linear-low} {\rm (ii)} as follows. 
\begin{eqnarray*}
I_1 &\leq &C\displaystyle\int_{0}^{\frac{t}{2}}(1+t-\tau)^{-\frac{n}{2}-\frac{1}{2}}\Big\|
\begin{pmatrix}
0 \\
m \otimes m
\end{pmatrix}
\Big\|_{L^1}d\tau\nonumber\\
&\leq & 
C\displaystyle\int_{0}^{\frac{t}{2}}(1+t-\tau)^{-\frac{n}{2}-\frac{1}{2}}(1+\tau)^{-\frac{n}{4}\cdot 2}d\tau M^2(t)\nonumber\\
&\leq &C(1+t)^{-\frac{n}{2}-\frac{1}{2}}M^2(t)\displaystyle\int_{0}^{\frac{t}{2}}(1+\tau)^{-\frac{n}{2} }d\tau \nonumber\\
&\leq &C(1+t)^{-\frac{n}{2}-\frac{1}{2}}\delta_1(t)M^2(t).
\end{eqnarray*}
$I_2$ can be estimated directly. Then we see from the estimates $I_1$ and $I_2$ that 
\begin{eqnarray}
\displaystyle\int_{0}^t\Big\| E_{1,m}(t,\tau)P_1 
\begin{pmatrix}
0 \\
\div(m \otimes m)
\end{pmatrix}
\Big\|_{L^\infty}d\tau\leq C(1+t)^{-\frac{n}{2}-\frac{1}{2}}\delta_1(t)M^2(t).\label{est-convection-higher}
 \end{eqnarray}
 Since other nonlinear term can be estimated similarly to \eqref{est-convection-higher}, 
 we have \eqref{m-partest}. 
\eqref{m-partest} together with \eqref{est-linear-lowinfty}, \eqref{est-linear-lowinfty2}, \eqref{est-linear-lowinfty3} and \eqref{nonlinearityest-low1}
imply that 
\begin{eqnarray}
\|\phi_1(t)\|_{L^\infty}
&\leq & C(1+t)^{-\frac{3n-1}{4}}\|u_{01}\|_{L^1}+C(1+t)^{q(n)}M^2(t),\label{estnonlinearlowketuron1}\\
\| m_1(t)\|_{L^\infty}&\leq & C(1+t)^{-\frac{n}{2}}\|u_{01}\|_{L^1}+C(1+t)^{-\frac{n}{2}-\frac{1}{2}}\delta_1(t)M^2(t).\label{estnonlinearlowketuron2}
\end{eqnarray}
In addition, we also obtain the following $L^2$ type estimate similarly based on Theorem \ref{est-linear-low} {\rm (iv)}.
\begin{eqnarray}
\|\nabla^k u_1(t)\|_{L^2}
&\leq & C(1+t)^{-\frac{n}{4}-\frac{k}{2}}\|u_{01}\|_{L^1}+C(1+t)^{-\frac{n}{4}-\frac{k}{2}}M^2(t)\label{estnonlinearlowketuron3}
\end{eqnarray}
for $k=0,1$. 

Concerning the $L^\infty$ type estimate for $u_\infty$, since the estimate in Theorem \ref{est-linear-high} has singularity at $t=0$ and the estimate of derivative of solutions has stronger singularity, we use $L^2$ energy estimate instead of using Theorem \ref{est-linear-high}. By Lemma \ref{lem2.1.} we see that 
$$
\|u_\infty(t)\|_{L^\infty} \leq C\|u_\infty\|_{H^{s+1}\times H^s}.
$$
Hence it is enough to obtain $L^2$ type estimate of $u_\infty$.  
As for the $L^2$ type estimate of $u_\infty$, using $L^2$ energy estimate stated in \cite[Proposition 5.4]{Tsuda-CNSK}, the following proposition is obtained for the high frequency part. 

\vspace{2ex} 

\begin{prop}\label{energyest}
Let $s$ be a nonnegative integer satisfying  $s\geq [n/2]+1$. Assume that 
\begin{eqnarray*}
u_{0\infty}=\trans(\phi_{0\infty},m_{0\infty})\in H^{s+1}_{(\infty)}\times H^s_{(\infty)},\\
F_\infty =\trans(F^1_\infty, F^2_\infty)\in L^2(0,T'; H^{s}_{(\infty)}\times H^{s-1}_{(\infty)})
\end{eqnarray*}
for all $T'>0$. Assume also that $u_{\infty}=\trans(\phi_{\infty},m_{\infty})$ satisfies
 \begin{eqnarray}
\left\{
\begin{array}{ll}
\partial_{t}u_{\infty}+Au_{\infty}=F_{\infty},\label{highparttimeeq}\\
u_{\infty}|_{t=0}=u_{0\infty}
\end{array}
\right.
\end{eqnarray}
and 
\begin{eqnarray*}
\phi_{\infty}\in C([0,T'];H^{s+1}_{(\infty)})\cap L^2(0,T';H^{s+2}_{(\infty)}),\ 
m_{\infty}\in C([0,T'];H^{s}_{(\infty)})\cap L^{2}(0,T';H^{s+1}_{(\infty)})
\end{eqnarray*}
for all $T'>0$. 
Then there exists an energy functional ${\cal E}[u_{\infty}]$ 
such that 
there holds the estimate 
\begin{eqnarray}\label{energy}
\frac{d}{dt}{\cal E}[u_{\infty}](t)
+{d}(\|\nabla \phi_{\infty}(t)\|_{H^{s+1}}^{2}+\|\nabla m_{\infty}(t)\|_{H^{s}}^{2})
 \leq
C\|F_{\infty}(t)\|_{H^{s}\times H^{s-1}}^2
\end{eqnarray} 
on $(0,T')$ for all $T'>0$. 
Here $d$ is a positive constant; 
$C$ is a positive constant independent of  $T'$;  
${\cal E}[u_{\infty}]$ is equivalent to $\|u_{\infty}\|_{H^{s+1}\times H^s}^2$, i.e, 
$$
C^{-1}\|u_{\infty}\|_{H^{s+1}\times H^s}^2
\leq {\cal E}[u_{\infty}]
\leq C\|u_{\infty}\|_{H^{s+1}\times H^s}^2; 
$$
and ${\cal E}[u_{\infty}](t)$ is absolutely continuous in $t\in [0,T']$ for all $T'>0$.
\end{prop}

\vspace{2ex} 

To apply the energy method, concerning estimate of the nonlinearity $P_\infty F(u)$ in $H^{s}\times H^{s-1}$, the following estimate is verified by direct computation based on Lemmas \ref{lem2.1.}-\ref{lem2.3.}. 

\vspace{2ex} 
\begin{lem}\label{nonlinearest-high}
It holds that 
$$
\|P_{\infty}F(u)\|_{H^{s}\times H^{s-1}}\leq C(1+t)^{-\frac{n}{2}-1}M^2(t) 
+C(1+t)^{-\frac{n}{2}}M(t)\|\nabla u_\infty\|_{H^{s+1}\times H^s}.
$$
\end{lem}

\vspace{2ex} 
We set $D[u_\infty]=\|\nabla \phi_{\infty}(t)\|_{H^{s+1}}^{2}+\|\nabla m_{\infty}(t)\|_{H^{s}}^{2}$. 
By \eqref{energy} and Lemma \ref{nonlinearest-high}, 
it is obtained that there exists a positive constant $c_3$ such that  
\begin{eqnarray}
\lefteqn{{\cal E}[u_{\infty}](t)+d \displaystyle\int_{0}^{t}e^{-c_3(t-\tau)}D[u_\infty](\tau)d\tau}\nonumber\\
&&\quad \leq e^{-c_3 t}{\cal E}[u_{\infty}](0) \nonumber\\
&&\qquad  +CM^4 (t)\displaystyle\int_{0}^{t}
e^{-c_3(t-\tau)}(1+\tau)^{-n-2}d\tau+CM^2 (t)\displaystyle\int_{0}^{t}
e^{-c_3(t-\tau)}(1+\tau)^{-n}D[u_\infty](\tau)d\tau\nonumber\\
&&\quad \leq e^{-c_3 t}{\cal E}[u_{\infty}](0) \nonumber\\
&&\qquad  +CM^4 (t)(1+t)^{-n-2}+CM^2 (t)\displaystyle\int_{0}^{t}
e^{-c_3(t-\tau)}(1+\tau)^{-n}D[u_\infty](\tau)d\tau.\label{estenergy}
\end{eqnarray}
${\cal D}[u_\infty]$ and $\tilde{\cal E}[u_\infty]$ are defined by 
\begin{eqnarray*}
{\cal D}[u_\infty](t)&=&(1+t)^{\frac{3n-1}{2}}\displaystyle\int_{0}^{t}
e^{-c_3(t-\tau)}D[u_\infty](\tau)d\tau, \\
\tilde{\cal E}[u_\infty](t) &=&\sup_{0\leq \tau \leq t}(1+\tau)^{\frac{3n-1}{2}}{\cal E}[u_{\infty}](\tau)
\end{eqnarray*}
for $n=2,3$ and 
\begin{eqnarray*}
{\cal D}[u_\infty](t)&=&(1+t)^{n+1}\displaystyle\int_{0}^{t}
e^{-c_3(t-\tau)}D[u_\infty](\tau)d\tau, \\
\tilde{\cal E}[u_\infty](t) &=&\sup_{0\leq \tau \leq t}(1+\tau)^{n+1}{\cal E}[u_{\infty}](\tau)
\end{eqnarray*}
for $n \geq 4$. 
We see from \eqref{estenergy} that 
\begin{eqnarray}
\tilde{\cal E}[u_\infty](t)+d{\cal D}[u_\infty](t) \leq C(\tilde{\cal E}[u_\infty](0)+M^4(t)+CM^2(t){\cal D}[u_\infty](t)).\label{energy-est-keturon}
\end{eqnarray}
Combining \eqref{estnonlinearlowketuron1}, \eqref{estnonlinearlowketuron2}, \eqref{estnonlinearlowketuron3} with \eqref{energy-est-keturon}, it is derived that 
\begin{eqnarray}
M^2(t) + d{\cal D}[u_\infty](t)\leq C(\|u_0\|_{(H^{s+1}\times H^s) \cap L^1}+M^4(t)+CM^{2}(t){\cal D}[u_\infty](t)). \label{nolinearketuron0}
\end{eqnarray}
Especially, we get 
$$
M(0)\leq C_1 \|u_0\|_{(H^{s+1}\times H^s)\cap L^1}.
$$
Since $M(t)$ is continuous in $t$, there exists a time $t_1>0$ such that 
$$
M(t)\leq 2 C_1 \|u_0\|_{(H^{s+1}\times H^s) \cap L^1}
$$
for $t\in [0,t_1]$. This together with \eqref{nolinearketuron0} derives that if $\epsilon_1 $ in Theorem {\rm (ii)} is sufficient small we have that there exists a positive constant $C_2$ such that 
$$
M(t)\leq C_2  \ \ \mbox{uniformly  for all}  \ \ t. 
$$ 
Consequently, 
Theorem \ref{thm-nonlinear} {\rm (ii)} is verified.  

Theorem \ref{thm-nonlinear} {\rm (i)} directly follows from \eqref{estnonlinearlowketuron1} and \eqref{estnonlinearlowketuron2}

The proof of Theorem \ref{thm-nonlinear} {\rm (iii)} is given as follows. By \eqref{low} we derive that 
\begin{eqnarray*}
\| u_j(t)\|_{L^1}\leq \| E_j(t) u_{0j}\|_{L^1}+ \displaystyle\int_{0}^{t}\| E_j(t,\tau)P_j F(u_1+u_\infty)(\tau) \|_{L^1}d\tau
\end{eqnarray*}
for $j=1, \infty$. 
We set 
\begin{eqnarray*}
J_1&=&\displaystyle\int_{0}^{t-1}\|E_j(t,\tau)P_j F(u_1+u_\infty)(\tau) \|_{L^1}d\tau, \\
 J_2&=&\displaystyle\int_{t-1}^{t}\| E_j(t,\tau)P_j F(u_1+u_\infty)(\tau) \|_{L^1}d\tau. 
\end{eqnarray*}
Concerning the estimates of $J_1$ and $J_2$, we derive the following estimates by direct computations based on $L^1$ type estimates in Theorems \ref{est-linear-low}  and \ref{est-linear-high}. 
\begin{eqnarray}
J_1 &\leq &C\displaystyle\int_{0}^{t-1} (t-\tau)^{\frac{n-1}{4}}\|P_j F(u_1+u_\infty)(\tau)\|_{L^1}d\tau\nonumber\\
&\leq &C\displaystyle\int_{0}^{t-1} (t-\tau)^{\frac{n-1}{4}}
(1+\tau)^{-\frac{n}{2}-\frac{1}{2}}M^2(t)d\tau\nonumber\\
&\leq &Ct^{\frac{n-1}{4}}M^2(t),\label{estj-1}\\
J_2 &\leq &C\displaystyle\int_{t-1}^{t}(t-\tau)^{-\frac{1}{2}-\sigma_0}\|P_j F(u_1+u_\infty)(\tau)\|_{L^1}d\tau\nonumber\\
&\leq &C\displaystyle\int_{t-1}^{t} (t-\tau)^{-\frac{1}{2}-\sigma_0}
(1+\tau)^{-\frac{n}{2}-\frac{1}{2}}M^2(t)d\tau\nonumber\\
&\leq &Ct^{-\frac{n}{2}-\frac{1}{2}}M^2(t).\label{estj-2}
\end{eqnarray}
By \eqref{estj-1}, \eqref{estj-2}, Theorems \ref{est-linear-low}  and \ref{est-linear-high}, it holds that 
\begin{eqnarray*}
\| u_j(t)\|_{L^1}\leq Ct^{\frac{n-1}{4}}(\| u_0\|_{L^1}+M^2(t))
\end{eqnarray*}
for $j=1,\infty$. Since $u=u_1+u_\infty$ and $M(t)$ is bounded by $u_0$ in Theorem \ref{thm-nonlinear} {\rm (ii)}, we obtain Theorem \ref{thm-nonlinear} {\rm (iii)}. This completes the proof. $\hfill\square$

\vspace{2ex}
\noindent {\bf Acknowledgements.} 
The first author is partly supported by Grants-in-Aid for Scientific Research with the Grant number: 16H03945. 
The second author is partly supported by Grant-in-Aid for JSPS Fellows with the Grant number: A17J047780. 


\end{document}